\documentclass[11pt]{article}
\usepackage[normalem]{ulem}
\usepackage{color}
\usepackage{multirow}
\usepackage{amsfonts,amssymb,a4, graphicx, bm, makeidx}
\usepackage{amsmath}
\usepackage{soul}
\usepackage{mathrsfs}
\usepackage[all]{xy}
\usepackage[latin1]{inputenc}
\usepackage[dvips]{geometry}
\def\theequation{\thesection.\arabic{equation}}

\newcommand{\qed}{\hfill\rule{3mm}{3mm}}
\newtheorem{corollary}{Corollary}[section]
\newtheorem{theorem}[corollary]{Theorem}
\newtheorem{lemma}[corollary]{Lemma}
\newtheorem{proposition}[corollary]{Proposition}
\newtheorem{definition}[corollary]{Definition}

\makeatletter \@addtoreset{equation}{section} \makeatother
\definecolor{GreenYellow}{cmyk}{0.15,0,0.69,0}
\definecolor{Yellow}{cmyk}{0,0,1,0}
\definecolor{Goldenrod}{cmyk}{0,0.10,0.84,0}
\definecolor{Dandelion}{cmyk}{0,0.29,0.84,0}
\definecolor{Apricot}{cmyk}{0,0.32,0.52,0}
\definecolor{Peach}{cmyk}{0,0.50,0.70,0}
\definecolor{Melon}{cmyk}{0,0.46,0.50,0}
\definecolor{YellowOrange}{cmyk}{0,0.42,1,0}
\definecolor{Orange}{cmyk}{0,0.61,0.87,0}
\definecolor{BurntOrange}{cmyk}{0,0.51,1,0}
\definecolor{Bittersweet}{cmyk}{0,0.75,1,0.24}
\definecolor{RedOrange}{cmyk}{0,0.77,0.87,0}
\definecolor{Mahogany}{cmyk}{0,0.85,0.87,0.35}
\definecolor{Maroon}{cmyk}{0,0.87,0.68,0.32}
\definecolor{BrickRed}{cmyk}{0,0.89,0.94,0.28}
\definecolor{Red}{cmyk}{0,1,1,0}
\definecolor{OrangeRed}{cmyk}{0,1,0.50,0}
\definecolor{RubineRed}{cmyk}{0,1,0.13,0}
\definecolor{WildStrawberry}{cmyk}{0,0.96,0.39,0}
\definecolor{Salmon}{cmyk}{0,0.53,0.38,0}
\definecolor{CarnationPink}{cmyk}{0,0.63,0,0}
\definecolor{Magenta}{cmyk}{0,1,0,0}
\definecolor{VioletRed}{cmyk}{0,0.81,0,0}
\definecolor{Rhodamine}{cmyk}{0,0.82,0,0}
\definecolor{Mulberry}{cmyk}{0.34,0.90,0,0.02}
\definecolor{RedViolet}{cmyk}{0.07,0.90,0,0.34}
\definecolor{Fuchsia}{cmyk}{0.47,0.91,0,0.08}
\definecolor{Lavender}{cmyk}{0,0.48,0,0}
\definecolor{Thistle}{cmyk}{0.12,0.59,0,0}
\definecolor{Orchid}{cmyk}{0.32,0.64,0,0}
\definecolor{DarkOrchid}{cmyk}{0.40,0.80,0.20,0}
\definecolor{Purple}{cmyk}{0.45,0.86,0,0}
\definecolor{Plum}{cmyk}{0.50,1,0,0}
\definecolor{Violet}{cmyk}{0.79,0.88,0,0}
\definecolor{RoyalPurple}{cmyk}{0.75,0.90,0,0}
\definecolor{BlueViolet}{cmyk}{0.86,0.91,0,0.04}
\definecolor{Periwinkle}{cmyk}{0.57,0.55,0,0}
\definecolor{CadetBlue}{cmyk}{0.62,0.57,0.23,0}
\definecolor{CornflowerBlue}{cmyk}{0.65,0.13,0,0}
\definecolor{MidnightBlue}{cmyk}{0.98,0.13,0,0.43}
\definecolor{NavyBlue}{cmyk}{0.94,0.54,0,0}
\definecolor{RoyalBlue}{cmyk}{1,0.50,0,0}
\definecolor{Blue}{cmyk}{1,1,0,0}
\definecolor{Cerulean}{cmyk}{0.94,0.11,0,0}
\definecolor{Cyan}{cmyk}{1,0,0,0}
\definecolor{ProcessBlue}{cmyk}{0.96,0,0,0}
\definecolor{SkyBlue}{cmyk}{0.62,0,0.12,0}
\definecolor{Turquoise}{cmyk}{0.85,0,0.20,0}
\definecolor{TealBlue}{cmyk}{0.86,0,0.34,0.02}
\definecolor{Aquamarine}{cmyk}{0.82,0,0.30,0}
\definecolor{BlueGreen}{cmyk}{0.85,0,0.33,0}
\definecolor{Emerald}{cmyk}{1,0,0.50,0}
\definecolor{JungleGreen}{cmyk}{0.99,0,0.52,0}
\definecolor{SeaGreen}{cmyk}{0.69,0,0.50,0}
\definecolor{Green}{cmyk}{1,0,1,0}
\definecolor{ForestGreen}{cmyk}{0.91,0,0.88,0.12}
\definecolor{PineGreen}{cmyk}{0.92,0,0.59,0.25}
\definecolor{LimeGreen}{cmyk}{0.50,0,1,0}
\definecolor{YellowGreen}{cmyk}{0.44,0,0.74,0}
\definecolor{SpringGreen}{cmyk}{0.26,0,0.76,0}
\definecolor{OliveGreen}{cmyk}{0.64,0,0.95,0.40}
\definecolor{RawSienna}{cmyk}{0,0.72,1,0.45}
\definecolor{Sepia}{cmyk}{0,0.83,1,0.70}
\definecolor{Brown}{cmyk}{0,0.81,1,0.60}
\definecolor{Tan}{cmyk}{0.14,0.42,0.56,0}
\definecolor{Gray}{cmyk}{0,0,0,0.50}
\definecolor{Black}{cmyk}{0,0,0,1}
\definecolor{White}{cmyk}{0,0,0,0}

\begin{document}
\def\theequation{\thesection.\arabic{equation}}

\def\blu{\color{Blue}}
\def\mag{\color{Maroon}}
\def\red{\color{Red}}
\def\green{\color{ForestGreen}}
\def\prob{{\rm Prob}}

\def\reff#1{(\protect\ref{#1})}

\let\a=\alpha \let\b=\beta \let\ch=\chi \let\d=\delta \let\e=\varepsilon
\let\f=\varphi \let\g=\gamma \let\h=\eta    \let\k=\kappa \let\l=\lambda
\let\m=\mu \let\n=\nu \let\o=\omega    \let\p=\pi \let\ph=\varphi
\let\r=\rho \let\s=\sigma \let\t=\tau \let\th=\vartheta
\let\y=\upsilon \let\x=\xi \let\z=\zeta
\let\D=\Delta \let\F=\Phi \let\G=\Gamma \let\L=\Lambda \let\Th=\Theta
\let\O=\Omega \let\P=\Pi \let\Ps=\Psi \let\Si=\Sigma \let\X=\Xi
\let\Y=\Upsilon

\global\newcount\numsec\global\newcount\numfor
\gdef\profonditastruttura{\dp\strutbox}
\def\senondefinito#1{\expandafter\ifx\csname#1\endcsname\relax}
\def\SIA #1,#2,#3 {\senondefinito{#1#2}
\expandafter\xdef\csname #1#2\endcsname{#3} \else
\write16{???? il simbolo #2 e' gia' stato definito !!!!} \fi}
\def\etichetta(#1){(\veroparagrafo.\veraformula)
\SIA e,#1,(\veroparagrafo.\veraformula)
 \global\advance\numfor by 1
 \write16{ EQ \equ(#1) ha simbolo #1 }}
\def\etichettaa(#1){(A\veroparagrafo.\veraformula)
 \SIA e,#1,(A\veroparagrafo.\veraformula)
 \global\advance\numfor by 1\write16{ EQ \equ(#1) ha simbolo #1 }}
\def\BOZZA{\def\alato(##1){
 {\vtop to \profonditastruttura{\baselineskip
 \profonditastruttura\vss
 \rlap{\kern-\hsize\kern-1.2truecm{$\scriptstyle##1$}}}}}}
\def\alato(#1){}
\def\veroparagrafo{\number\numsec}\def\veraformula{\number\numfor}
\def\Eq(#1){\eqno{\etichetta(#1)\alato(#1)}}
\def\eq(#1){\etichetta(#1)\alato(#1)}
\def\Eqa(#1){\eqno{\etichettaa(#1)\alato(#1)}}
\def\eqa(#1){\etichettaa(#1)\alato(#1)}
\def\equ(#1){\senondefinito{e#1}$\clubsuit$#1\else\csname e#1\endcsname\fi}
\let\EQ=\Eq

\def\pp{{\bm p}}\def\pt{{\tilde{\bm p}}}


\def\\{\noindent}
\let\io=\infty
\def\ee{\end{equation}}
\def\be{\begin{equation}}
\def\bea{\begin{eqnarray}}
\def\eea{\end{eqnarray}}

\def\VU{{\mathbb{V}}}
\def\EE{{\mathbb{E}}}
\def\N{\mathbb{N}}
\def\U{\mathbb{U}}
\def\GI{{\mathbb{G}}}
\def\TT{{\mathbb{T}}}
\def\C{\mathbb{C}}
\def\CC{{\mathcal C}}
\def\KK{{\mathcal K}}
\def\II{{\mathcal I}}
\def\LL{{\cal L}}
\def\RR{{\cal R}}
\def\SS{{\cal S}}
\def\NN{{\cal N}}
\def\HH{{\cal H}}
\def\GG{{\cal G}}
\def\PP{{\cal P}}
\def\AA{{\cal A}}
\def\BB{{\cal B}}
\def\FF{{\cal F}}
\def\v{\vskip.1cm}
\def\vv{\vskip.2cm}
\def\gt{{\tilde\g}}
\def\E{{\mathcal E} }
\def\EI{{\mathbb E} }
\def\I{{\rm I}}
\def\rfp{R^{*}}
\def\rd{R^{^{_{\rm D}}}}
\def\ffp{\varphi^{*}}
\def\ffpt{\widetilde\varphi^{*}}
\def\fd{\varphi^{^{_{\rm D}}}}
\def\fdt{\widetilde\varphi^{^{_{\rm D}}}}
\def\pfp{\Pi^{*}}
\def\pd{\Pi^{^{_{\rm D}}}}
\def\pbfp{\Pi^{*}}
\def\fbfp{{\bm\varphi}^{*}}
\def\fbd{{\bm\varphi}^{^{_{\rm D}}}}
\def\rfpt{{\widetilde R}^{*}}
\def\A{{{\mathcal O}}}
\def\ef{\mathfrak{f}}
\def\Ti{\mathfrak{T}}
\def\Mi{\mathfrak{M}}
\def\mm{\mathrm{ \mathbf{m}}}

\def\tende#1{\vtop{\ialign{##\crcr\rightarrowfill\crcr
              \noalign{\kern-1pt\nointerlineskip}
              \hskip3.pt${\scriptstyle #1}$\hskip3.pt\crcr}}}
\def\otto{{\kern-1.truept\leftarrow\kern-5.truept\to\kern-1.truept}}
\def\arm{{}}
\font\bigfnt=cmbx10 scaled\magstep1

\newcommand{\card}[1]{\left|#1\right|}
\newcommand{\und}[1]{\underline{#1}}
\def\1{\rlap{\mbox{\small\rm 1}}\kern.15em 1}
\def\ind#1{\1_{\{#1\}}}
\def\bydef{:=}
\def\defby{=:}
\def\buildd#1#2{\mathrel{\mathop{\kern 0pt#1}\limits_{#2}}}
\def\card#1{\left|#1\right|}
\def\proof{\noindent{\bf Proof. }}
\def\qed{ \square}
\def\trp{\mathbb{T}}
\def\trt{\mathcal{T}}
\def\Z{\mathbb{Z}}
\def\be{\begin{equation}}
\def\ee{\end{equation}}
\def\bea{\begin{eqnarray}}
\def\eea{\end{eqnarray}}
\def\kk{{\bf k}}
\def\Ti{\mathfrak{T}}
\def\Mi{\mathfrak{M}}
\def\begn{\begin{aligned}}
\def\egn{\end{aligned}}
\def\ti{{\rm\bf  t}}\def\mi{{\rm\bf m}}
\def\Va{{V^a_{\rm h.c.}}}
\def\Re{{\mathbb{R}}}
\def\T{{\mathcal{T}}}
\def\hL{{\L}}
\def\ev{\mathfrak{e}}
\def\obj{{\rm supp}}\def\fa{\FF}
\def\E{{\cal E}}
\def\EA{{E_\A}}
\def\0{\emptyset}
\def\Ni{\overline{\N}}
\def\zt{{\tilde z}}

\title{On the independent set polynomial of  graphs and claw-free graphs
}

\author{
\\
Paula M. S. Fialho$^1$,  Aldo Procacci$^2$. \\
\\
\small{$^1$Departamento de Ci\^encia da Computa\c{c}\~ao UFMG, }
\small{30161-970 - Belo Horizonte - MG
Brazil}\\
\small{paulamsfialho@gmail.com}\\
\small{$^2$Departamento de Matem\'atica UFMG,}
\small{ 30161-970 - Belo Horizonte - MG
Brazil}\\
\small{aldo@ufmg.br}
}

\maketitle

\begin{abstract}
{ We present two new contributions to the study of the independence polynomial $Z_\GI(z)$ of a finite simple graph $\GI = (\VU,\EE)$.  First, we provide an improved lower bound for the zero-free region of $Z_\GI(z)$  for the important class of claw-free graphs. Our bound exceeds the classical Shearer radius and it is derived through a refined application of the Fern\'andez-Procacci criterion using properties of the local neighborhood structure in claw-free graphs.
Second, we establish a novel combinatorial expression for $Z_\GI(z)$, inspired by the connection with the abstract polymer gas models in statistical mechanics, which offers a new structural interpretation of the polynomial and may be of independent interest.
These results strengthen the connection between statistical physics, combinatorics, and graph theory, and suggest new approaches for analytic exploration.}
\end{abstract}

\numsec=2\numfor=1

\section{Introduction}

In this note $\GI=(\VU, \EE)$  will denote a simple graph with vertex set $\VU$ and  edge set $\EE$.
Specifically, $\VU$ is a finite set  and $\EE\subset [\VU]^2$ where $ [\VU]^2$
is the set of
all subsets of $\VU$ with cardinality 2.
We recall that $R \subseteq \VU$
is an {\it independent set of $\GI$} if for all $\{x,y\}\subset R$ we have that
 $\{x,y\}\notin \EE$ and we denote by $I(\GI)$ the set of all independent sets of $\GI$.

\\Given $z\in \C$, the {\it univariate independent set polynomial} of $\GI$ is defined by
\be\label{indset}
Z_\GI(z)=\sum_{S\in I(\GI)}z^{|S|}
\ee
and, given $\bm z=\{z_v\}_{v\in \VU} \in \mathbb{C}^{|\VU|}$, its multivariate generalization is given by
\begin{equation}\label{eq:zmulti}
Z_\GI(\bm z)= \sum_{S\in I(\GI)}\prod_{v\in S}z_v.
\end{equation}

\\This polynomial,  also referred to simply as the independence polynomial, is an important graph polynomial that arises in many contexts in combinatorics, mathematical physics and computer science. Indeed, in the context of mathematical physics,
the multivariate version $Z_\GI(\bm z)$ of the independence polynomial given in \reff{eq:zmulti} is also known as the grand-canonical partition function of the abstract polymer gas. From this perspective,  the vertices of $\GI$ are interpreted as ``polymers'',  where each polymer
$v\in \VU$ is assigned an activity $z_v \in \mathbb{C}$, and two polymers $v,v' \in \VU$ are {\it incompatible} if either $v=v'$ or $\{v,v'\}\in \EE$.

\\The abstract polymer gas, originally introduced in \cite{KP}, is a fundamental model in statistical mechanics. Its importance is due to the fact that   the partition function of the vast majority of spin systems  in  a lattice can be rewritten in terms of the partition function of an abstract polymer gas  (this can practically always be done in the high temperature regime and often in the low temperature regime).
Thus, investigating the zeros of
$Z_\GI(\bm z)$  is a crucial issue, since the existence of phase transitions in any statistical mechanics model is directly related to the zero-free region of the model's partition function.
The best known bounds for the zero-free region of the multivariate independence polynomial, were established in \cite{FP}.
Given a vertex $v \in \VU$,{  let $\G_\GI(v)$ denote its neighborhood in $\GI$, defined as:
$\G_\GI(v)= \{ v' \in \VU:~ \{v, v'\} \in \EE \}$}. Then the following theorem holds (see \cite{FP}).

\begin{theorem}[Fern\'andez-Procacci]\label{FPC}
Let $\bm\mu= \{\mu_v\}_{v\in \VU}$ be a collection of nonnegative numbers and let
$\bm r^*= \{\\r^*_v\}_{v\in \VU}$ such that
 \be\label{fp}
  r^*_v\equiv \frac{\mu_v}{\m_v+\ph_v(\bm \mu)},~~~~~~~~~~\forall v\in \VU
 \ee
 with
\be\label{FPfu}
\ph_v(\bm \mu)= \;\sum_{S\subset \G_\GI(v)\atop S\in I(\GI)}\prod_{v'\in S}\m_{v'}
\ee
then     $Z_\GI(\bm z)\neq 0$ in the polydisc $|\bm z|\le \bm r^*$ and therefore $Z_\GI (-\bm r^*)>0$.
\end{theorem}

\\Similarly, the univariate independence polynomial $Z_\GI(z)$ given in (\ref{indset}) can  be interpreted as  the grand canonical partition function of the self-repulsive hard core lattice gas. In this model identical particles with activity $z$ occupy the vertices $\VU$ of $\GI$, each vertex can hold at most one particle (self-repulsion), and no two particles can occupy adjacent vertices (hard-core constraint).
There are several results concerning the location of the zeros of $Z_\GI(z)$. In particular, Shearer \cite{Sh} (see also \cite{SS}) proved  that the univariate independence polynomial  $Z_\GI(z)$ of a graph $\GI$ with maximum degree $\D$ is free of zeros inside  the complex disk centered at the origin with radius
\be\label{shea}
r_\D={(\D-1)^{\D-1}\over \D^\D}.
\ee
The radius of the  Shearer disk is optimal in the sense that   if $\GI$  is the complete rooted tree
with branching factor $\D-1$ and depth $n$, then $Z_\GI(z)$  has a negative
real zero that approaches  $r_\D$ from the left as $n\to \infty$. The zero-free region in the complex plane of the univariate independence
polynomial has been  further investigated in recent works  including  \cite{BBP}, \cite{BCSV}, and \cite{PeR}. These studies have progressively expanded  the zero-free region  primarily towards the right half of the complex  plane.
{Nevertheless, the Shearer radius
$r_\D$ remains the best known lower bound for the minimal absolute value of the first zero.}


\\The independence polynomial of a graph also plays an important role in the probabilistic method in combinatorics, due to the connection --pointed out by Scott and Sokal \cite{SS}--
between the  abstract polymer gas and the famous Lov\'asz local lemma (LLL) originally  introduced by Erd\H os and Lov\'asz in \cite{EL}.
The LLL provides a sufficient condition on the probabilities $\mathbf{p}=\{p_\ev\}_{\ev\in \FF}$ of a finite family $\FF$ of undesirable (bad) events in some probability space, ensuring that none event in $\FF$ occurs. This condition is formulated in terms of a graph $\GI$
with vertex set $\FF$ called the dependence graph associated to the family $\FF$  of bad events (see e.g. \cite{BFPS}, \cite{SS} and references therein). As shown by Scott and Sokal, the necessary and sufficient condition  given by Shearer in \cite{Sh} for the thesis of the LLL to hold is equivalent to requiring  that the independence polynomial $Z_\GI(z)$ has no zeros in the  poly-disk $|\bm z|\le \mathbf{p}$ in $\C^{|\FF|}$.


\\Finally, we highlight that
the location of the zeros of  $Z_\GI(z)$ has become recently relevant also in theoretical computer science. Indeed, zero-free regions of the independence polynomial are crucial in Barvinok's method for deterministic approximate counting \cite{B}, as well as in the design of deterministic algorithms via interpolation and complex analysis \cite{PR}.


\\In the present paper we provide two observations on the independence polynomial of a graph $\GI$ that we believe may provide new structural and analytic insights, potentially leading to further developments in both theory and applications.
In Section \ref{Three}, we establish  (in Theorem \ref{clawteo}) a new bound that extends beyond the Shearer radius of the zero-free region of $Z_\GI(z)$ when $\GI$ is a claw-free graph.
The relevance and motivation behind studying this  class  of graphs will be discussed there.
In Section \ref{two}, we derive (in Theorem \ref{tea}) a novel formulation of $Z_\GI(z)$, based on its connection with statistical mechanics and inspired by the recent work \cite{FJP}.



\section{Zero-free regions for the independence polynomial of claw-free graphs}\label{Three}

\\Let us start this section by noting  that from  Theorem \ref{FPC}, the following corollary regarding the univariate independence polynomial follows immediately.

\begin{corollary}\label{corofp}
The univariate independence polynomial on a graph $\GI=(\VU,\EE)$ is free of zeros for any $\m>0$ such that
\be\label{fpw}
|z|\le\min_{v\in \VU}{\m\over \mu+\ph_v(\mu)}
\ee
where
\be\label{fist}
\ph_v(\mu)=\sum_{S\subset \G_\GI(v)\atop S\in I(\GI)}\m^{|S|}.
\ee
\end{corollary}
\vv
Observe that when $\GI$ is  a triangle-free graph, i.e. such that for any
$v\in \VU$ the neighborhood  $\G_\GI(v)$ is an independent set   (e.g.  a tree or a bipartite graph)
the function  \reff{fist} in Corollary \ref{corofp}
becomes
$$
\ph_v(\mu)=\sum_{S\subset \G_\GI(v)}\m^{|S|}= (1+\mu)^{d_v}
$$
where $d_v$ is the degree of the vertex $v$. Therefore if $\GI$ is a triangle free  graph with maximum degree $\D$, we get from Corollary \ref{corofp} that $Z_\GI(z)$ is free of zeros as soon as
$$
|z|\le {\mathfrak{r}_\D}:= \max_{\m>0}{\mu\over \mu+(1+\mu)^\D}={(\D-1)^{\D-1}\over \D^\D}\cdot {1\over 1+ {(\D-1)^{\D-1}\over \D^\D}}
$$
The  radius ${\mathfrak{r}_\D}$ of the zero-free disk for $Z_\GI(z)$  deduced from Corollary \ref{corofp} is slightly smaller than the  Shearer radius  $ r_\D$ given in \reff{shea}, and ${\mathfrak{r}_\D} \to r_\D$ as $\D\to \infty$. Nevertheless,  the function $\ph_v(\mu)$ defined in \reff{fist}, which depends strongly on the structure of the vertex's neighborhood,  suggests that Corollary \ref{fpw} may be used  to
improve on the Shearer bound  when $\GI$ is a  non triangle-free graph.

\vv

\\An important and widely studied class of non-triangle-free graphs is the class of claw-free graphs.
A graph $\GI=(\VU, \EE)$ is  claw-free  if $\GI$ does not contain the complete bipartite graph $K_{1,3}$ as an induced subgraph.
\\Claw-free graphs  have been intensively studied (see e.g. \cite{PaR}, \cite{Mi}, \cite{CS},  the surveys  \cite{FFR}, \cite{CSbook}, and also the sequence of papers \cite{CS1}-\cite{CS7}), with the initial motivation that any line graph\footnote{
Given a graph $G=(V,E)$, the line graph $L(G)$ of $G$
is the graph with vertex set $E$ and a pair  $\{e,e'\}$ belongs to the edge set of $L(G)$ if and only if they are adjacent (share a vertex) in $G$. The class of all line graphs is a proper subclass of claw-free graphs.} is claw-free.
Actually,  if $\GI$ is
a line graph, then a famous result
by Heilmann and Lieb  on the so-called Monomer-dimer model \cite{HL} implies that all roots of $Z_\GI(z)$ are real negative. Moreover, letting $\l_1(\GI)$ denote the real root of $Z_\GI(z)$  closest to the origin, the authors proved in \cite{HL} that for any line graph $\GI$ with maximum degree $\D$,
\be\label{Hlib}
\l_1(\GI)\le -{1\over 2\D},
\ee
whose absolute value is larger than the Shearer radius ($r_\D\approx {1\over e\D}$ as $\D\to\infty$).

\\In 2007, Chudnovsky and Seymour extended the Heilmann and Lieb results by proving, in their celebrated paper \cite{CS}, that the roots of the independence polynomial of any claw-free graph are all real (and hence negative). However, in \cite{CS} the authors do not provide a similar bound to \reff{Hlib} for $\l_1(\GI)$, the real root closest to the origin of the independence polynomial of a claw-free graph $\GI$ with maximum degree $\D$.

\\With the aim of obtaining a bound {\it \`a la} Heilmann and Lieb (i.e. bound \reff{Hlib}) for all claw-free graphs, we need to investigate the properties of the function
$\ph_v(\mu)$. The structure of the neighborhood $\G_\GI(v)$ of a vertex $v$ in a claw-free graph $\GI$ is such that if $S\in \G_\GI(v)$ is independent, then necessarily  $|S|\le 2$. Therefore, denoting by $s_v$ the number of independent pairs in $\G_\GI(v)$ and by
$\Phi=\max_{v\in \VU}s_v$,  we have that in a claw-free graph
\be\label{gene}
\ph_v(\mu)\le 1+(\D+1)\m +\Phi\m^2.
\ee

\\Bound \reff{gene} together with  Corollary \ref{corofp}  immediately implies the following
theorem providing a general upper bound for  $\l_1(\GI)$ when $\GI$ is claw free.

\vv
\begin{theorem}\label{clawteo}
Let $\GI=(\VU,\EE)$ be a claw-free graph with maximum degree $\D$ and  let  $\Phi$ be the maximum number of independent pairs in the  neighborhood of any vertex of $\GI$. Then the univariate independence polynomial $Z_\GI(z)$ of $\GI$ is free of zeros as soon as
\be\label{cfg}
|z|\le {1\over \D+1+2\sqrt{\Phi}}
\ee
and therefore
\be\label{bdphi}
\l_1(\GI)\le - {1\over \D+1+2\sqrt{\Phi}}
\ee
\end{theorem}

\vv

\vv
\\{\bf Proof}.  By Corollary \reff{corofp} and bound \reff{gene}, $Z_\GI(z)$ is free of zeros for any $\m>0$ such that
$$
|z|\le {\mu\over 1+(\D+1)\m+\Phi\mu^2}.
$$
and optimizing in $\m$ we get \reff{cfg}.

~~~~~~~~~~~~~~~~~~~~~~~~~~~~~~~~~~~~~~~~~~~~~~~~~~~~~~~~~~~~~~~~~~~~~~~~~~~~~~~~~~~~~~~~~~~~~~~~~~~~~~~~~~~~~~~~~~~~~~~~~~~~~~~~~~~~~~~~~~~~~~~~~~~~~~~~~$\Box$

\\In order to compare \reff{bdphi} with the Heilmann-Lieb bound \reff{Hlib}, we also  provide an alternative (and in general worst) version of bound \reff{bdphi} in terms of the maximum degree $\D$ only. Mantel's theorem \cite{Ma} states that any triangle-free graph with $n$ vertices contains at most ${n^2\over 4}$ edges. Thus, as the complement of a triangle-free graph is a claw-free graph,   the number of independent pairs in the neighborhood $\G_\GI(u)$ of any vertex $u$ of the claw-free graph $\GI$ with maximum degree $\D$ cannot exceed
  ${\D^2\over 4}$. In other words, we have that $\Phi\le {\D^2\over 4}$. Then, Theorem \ref{clawteo} immediately implies the following  corollary providing
a  bound for   $\l_1(\GI)$
that depends only on $\D$ and asymptotically matches the Heilmann-Lieb bound.

\begin{corollary}\label{clawcoro}
Let $\GI=(\VU,\EE)$ be a claw-free graph with maximum degree $\D$. Then
\be\label{bdelta}
\l_1(\GI)\le - {1\over 2\D+1}.
\ee
\end{corollary}

\vv

\vv

\\We now present, as an example that will be useful later, the so-called Sch\"afli graph.  The    Sch\"afli graph $G_{\rm sch}$ is a $16-$regular claw-free graph (see 10.10 in \cite{Br} and \cite{LR}). Direct calculations show that the smallest (in modulous) root of the independence polynomial of $G_{\rm sch}$ is \(\lambda_1(G_{\rm sch})=-0.048706 \). The neighborhood of any vertex of  the Sch\"afli graph  is isomorphic to the complement of the
 Clebsch graph, which is a 10-regular graph with  16 vertices and thus 80 edges (see 10.7 in \cite{Br}). Therefore,  for any vertex  $v$ of $G_{\rm sch}$ we have that
$\ph_v(\mu)={1+17\m+40\mu^2}$ and then
 Theorem \ref{clawteo} gives \(\lambda_1(G_{\rm sch}) \le-0.03373\), while Corollary \ref{clawcoro} gives the (slightly worse) bound
\(\lambda_1(G_{\rm sch}) \le-{1\over 33}\).

\\Theorem \ref{clawteo} and Corollary \ref{clawcoro} are not the first results on the smallest root \( \l_1(\GI) \) of the  univariate
independence polynomial of a claw-free graph $\GI$. Indeed, in a recent paper  \cite{LR} Leake and Ryder  investigate the independence polynomial from the perspective of multivariate stability theory and  also  derive the following upper bounds on  \( \l_1(\GI) \).

\begin{proposition}[Proposition 5.3 of \cite{LR}]\label{proplery} Given any claw-free graph $\GI$, we have
\be\label{lery}
\l_1(\GI) < -\frac{1}{4 \cdot \max\{\o - 1, \d\}},
\ee
where \( \o \) is the clique number of \( \GI \) and \( \d \) is the minimum degree.
\end{proposition}

\\The authors also show a tighter bound for a subclass of the claw-free graphs (i.e. those containing a simplicial clique) which includes all line graphs. They refer to  graphs in  this class as  ``simplicial graphs" (see \cite{LR}).


\begin{proposition}[Corollary 4.9 of \cite{LR}]\label{corolery2}
Given a simplicial graph $\GI$, we have that
\be\label{lery2}
\lambda_1(\GI) \le -{1\over 4(\omega - 1)}.
\ee
\end{proposition}
\vv

\\Observe that, if \( \GI \) is the line graph of a graph \( H \) with maximum degree \( d \), then \( \omega = d \), and the maximum degree of $\GI$ is  \( \Delta = 2(d - 1) = 2(\omega - 1) \) and thus their bound \reff{lery2} can be rewritten as \( \lambda_1(\GI) < -1/2\Delta \), which coincides with the Heilmann-Lieb bound.

\\However, their general bound  \reff{lery} valid for all claw-free graph is not directly comparable with ours, as  \reff{lery} depends on  the clique number and the minimum degree, whereas  \reff{bdphi} is based on the maximum degree and the maximum number of independent pairs in a neighbor.  The
relationship between these four parameters in a claw-free graph is not at all immediate. However,  we would like to stress that our $\D$-dependent bound \reff{bdelta} is worse than \reff{lery} only if $\o\le {\D\over 2}+1$  and $\d\le {\D\over 2}$.
\vv

\\Leake and Ryder also observe that  bound \reff{lery2}, valid for all simplicial graphs, may fail for non simplicial  claw-free graphs  and they give as an example the above mentioned Sch\"afli graph. Applying Proposition \ref{proplery} on the Sch\"afli graph $G_{\rm sch}$, which is such that  $\D=\d=16$ and $\o=6$ (see 10.10 in \cite{Br}), we get
 \(\lambda_1(G_{\rm sch}) \le -1/64=-0.015625\). However, Sch\"afli graph is not simplicial, and bound (\ref{lery2}) in Proposition \ref{corolery2}  would (wrongly) produce $\lambda_1(G_{\rm sch})\le -{1\over 4(\o-1)}=-0.05$.

\begin{table}[h!]
\centering
\renewcommand{\arraystretch}{1.2}
\begin{tabular}{|c|l|c|}
\hline
\textbf{Source} & \textbf{Technique} & \textbf{Bound / Value} \\
\hline
--           & Exact value (computed)           & $\lambda_1(G_{\rm sch})={-0.0487057}$ \\
\hline
FP           & Theorem~\ref{clawteo}            & $\lambda_1(G_{\rm sch}) \leq -0.03373$ \\
FP           & Corollary~\ref{clawcoro}         & $\lambda_1(G_{\rm sch}) \leq -\frac{1}{33} \approx -0.0303$ \\
\hline
LR           & Proposition \ref{proplery}             & $\lambda_1(G_{\rm sch}) \leq -\frac{1}{64} = -0.015625$ \\
\hline
\end{tabular}
\caption{Comparison between the exact value, Leake and Ryder bounds (LR), and our results (FP) for the Schl\"afli graph.}
\label{tab:schlafli-bounds}
\end{table}


\section{An alternative representation of $Z_\GI(z)$}\label{sec2}\label{two}

{In this section we  provide a novel combinatorial expression for the independence polynomial $Z_G(z)$, inspired by its connections with statistical mechanics, as exposed in introduction.}
We begin by establishing some notations  and recalling some definitions regarding graphs, trees and forests.
Given $n\in \N$, we set $[n]=\{1,2,\dots,n\}$ and
if  $U$ is a finite set, $|U|$ denotes its  cardinality and  $\mathcal{P}(U)$ denotes  the set of all subsets of $U$. Hereafter the symbol $\biguplus$ will denote the disjoint union.

\\A graph $\GI=(\VU, \EE)$
is {\it connected}
if for any pair $B, C$ of  subsets of $\VU$ such that
$B\biguplus C =\VU$, there is an edge  $e\in \EE$ such
that $e\cap B\neq\emptyset$ and $e\cap C\neq\emptyset$.
A subgraph of $\GI$ is a  graph $G=(R,E)$  where  $R \subseteq \VU$ and $E \subseteq \EE$ such that $\{x,y\}\subset R$ for any $\{x,y\}\in E$.
A connected component of $\GI$  is a maximal connected subgraph of $\GI$.
Given a set $R\subset \VU$ we denote by $\GI|_{R}$ the (induced) subgraph of $\GI$ with vertex set $R$
and edge set $\EE|_R=\{\{x,y\}\in \EE: \{x,y\}\subset R\}$.
A subset $R\subset \VU$ is said to be connected if $\GI|_{R}$ is connected; we denote by $\mathrm{C}_\VU$ the set of all connected subsets of $\GI$ with cardinality greater than one, i.e.
$$
\mathrm{C}_\VU=\{R\subset \VU:\; \GI|_R \; \mbox{is connected and} \;|R|\ge 2\}.
$$

\\Given {$E\subset \EE$, we set
$V_E=\{x\in \VU: x\in e ~{\rm for ~some~} e\in E\}$ and $\GI|_E$ denotes the subgraph of $\GI$ with vertex set $V_E$ and edge set $E$.
A non-empty subset $E\subset \EE$ is connected if the graph $\GI|_E$ is connected.
Any  non-empty $E\subset \EE$ can be written in a unique way, for some integer $k\in \N$,  as $E=\biguplus_{i=1}^kE_i$ in such a way that  $E_i$
is connected for all $i\in [k]$ and $V_{E_i}\cap V_{E_j}=\0$ for all $\{i,j\}\subset [k]$; the subsets $E_1, \dots, E_k$ are called the connected components of $E$.
A connected set $\t\subset \EE$ such that $|V_\t|=|\t|+1$ is called a tree in $\GI$.
{Given $E\subset \EE$,} a subset $E'\subset E$ is called a spanning subset of $E$ if
$V_{E'}=V_{E}$.
}
\\Given  $R\in \mathrm{C}_\VU$,  we denote by  $\mathcal{C}_{R}$ (resp. $\mathcal{T}_{R}$)  the set of  all spanning connected subsets (resp. spanning trees) of $\EE|_R$.

{
\\A forest in $\GI=(\VU,\EE)$ is {either the empty set or} a {non-empty} subset $F\subset \EE$ whose  connected components are trees.
We denote by $\mathfrak{F}_\GI$ the set of all forests in $\GI$.
Any  non-empty forest $F\in \mathfrak{F}_\GI$ can be written as  $F=\biguplus_{i=1}^k\t_i$, for some $k\in \N$,  with  $\t_1,\dots, \t_k$ non-trivial trees such that
$V_{\t_i}\cap V_{\t_j}=\0$ for all pairs $\{i,j\}\subset [k]$.
We denote  by $\|F\|$ the number of trees  forming $F$. Please note that   $V_F=\biguplus_{i=1}^kV_{\t_i}$ and  $|V_F|=|F|+\|F\|$.
Observe that a tree $\t\subset \EE$ is just a non-empty forest constituted by  a single connected component. We denote by $\mathfrak{T}_\GI$ the set of all forest in $\GI$ which are single trees (i.e.  $\mathfrak{T}_\GI=\{F\in \mathfrak{F}_\GI:~\|F\|=1\}$).}


\\Given $R\in \mathrm{C}_\VU$, {and} $E, E'\in \mathcal{C}_R$ such that $E \subset E'$, let us define the Boolean interval
$$[E, E']=\{E''\in \mathcal{C}_R:~ E \subset E'' \subset E'\}.$$

\\Let us now introduce the notion of a {\it partition scheme in $\GI$}, { a systematic way to partition the set of all  forests into non-overlapping classes.}

{\begin{definition}\label{defia}
A partition scheme in $\GI$ is a map  $\mm: \mathfrak{F}_\GI\to \mathcal{P}(\EE)$
 such that:
\vv
(1) $\mm(\0)=\0$
\vv
(2) $\t\subset \mm(\t)$ and $V_\t=V_{\mm(\t)}$  for each $\t\in {\mathfrak{T}}_\GI$;
\vv
(3) $\mathcal{C}_R=\biguplus_{\tau\in \mathcal{T}_R}[\tau, \mm(\tau)]$ for each $R\in \mathrm{C}_\VU$;
\vv
(4) $\mm(F)= \bigcup_{i=1}^{\|F\|}\mm(\t_i)$ for any non-empty forest  $F\in \mathfrak{F}_\GI$ such that $F={\biguplus_{i=1}^{ \|F\|}\t_i}$.
\end{definition}}

\\Given any partition scheme $\mi$ in $\GI$ and given $R\in \mathrm{C}_\VU$, we set
$\mathcal{T}_{R}^{\mi}=\{ \t \in \mathcal{T}_{R}:~ \mi(\t)=\t\}$  and

\begin{equation}\label{def:fmi}
\FF_{\GI,\mi}=\{F={\biguplus_{i=1}^{ \|F\|}\t_i} \in \mathfrak{F}_\GI:~ \t_i \in \mathcal{T}_{R}^{\mi}~\text{for all} ~ i\in [k] \}.
\end{equation}
Several partition schemes are available; see e.g.  \cite{pen67}, \cite{FP}, \cite{JPS} for the so-called Penrose Partition scheme, and \cite{W}, \cite{PY}, \cite{FJP} for the ``minimal-tree" partition scheme. See also \cite{GS} and    Section 4 of \cite{SS}  for other examples of partition schemes.


\vv

\\Now we have all the elements to state our result.

\begin{theorem}\label{tea}
Consider a finite simple graph $\GI=(\VU, \EE)$. Let $Z_{\GI}(z)$ represent its independent set polynomial, with $z \in \mathbb{C}$ and let $\mi$ be any partition scheme in $\GI$. Then the following identity holds
\be\label{teoz}
Z_\GI(z) =(1+z)^{|\VU|}\sum_{F\in \FF_{\GI,\mi} }(-1)^{|F|}\left({z\over 1+z}\right)^{|V_F|}.
\ee
\end{theorem}
\vv

%
\\The proof of theorem \ref{tea} will be given in the next section.
Possible applications, such as, e.g. further improvements of the zero free region
of the independence polynomial
 of specific classes of graphs based on the recent approaches to zero-freeness of the chromatic polynomial given in \cite{JPR,FJP2,BR},
 will be investigated in a separated paper.
Here we just comment that the alternative expression of the independence polynomial given by the r.h.s. of \reff{teoz} could be useful to compute the independence polynomial at certain values of of $z$. In particular, the important value of the independence polynomial of a graph $\GI$ at $z=-1$, which is is related to the reduced Euler charateristic of the independence complex of $\GI$
(see e.g. \cite{LM} and reference therein), according to Theorem \ref{tea}, is  equal to
\be\label{z-1}
Z_\GI(z=-1)=\sum_{F\in \FF_{\GI,\mi}\atop V_F=|\VU| }(-1)^{\|F\|}.
\ee
I.e., $Z_\GI(z=-1)$ is  the number of $\mi$-invariant (where $\mi$ is any partition scheme
in $\GI$) {spanning} forests\footnote{A forest in $\GI$ is spanning if it does not contain isolated vertices} of $\GI$ composed by an even number of trees minus the number of $\mi$-invariant spanning forests of $\GI$ formed by an odd number of trees.

\\Theorem \ref{tea} also gives interesting alternative expressions for $Z_\GI(z=-1/2)$ and $Z_\GI(z=1)$ in terms of $\mi$-invariant forests. Namely,
$$
Z_\GI(z=-1/2)={1\over 2^{|\VU|}}\sum_{F\in \FF_{\GI,\mi}}(-1)^{\|F\|}
$$
and
$$
Z_\GI(z=1)={2^{|\VU|}}\sum_{F\in \FF_{\GI,\mi}}(-1)^{\|F\|}.
$$

\subsection{ Proof of Theorem \ref{tea}}
\vv

We define a pair potential from $V: [\VU]^2 \to \{0,+\infty\}$ in such a way that
\begin{equation}\label{eq:V}
V(\{x,y\})=\begin{cases} 0 &{\rm if} ~\{x,y\}\not\in\EE\\
+\infty  &{\rm if} ~\{x,y\}\in\EE
\end{cases}
\end{equation}
so that the independence polynomial  (\ref{indset}) can be written as
$$
Z_\GI(z)=\sum_{S\subseteq \VU}z^{|S|}e^{-\sum_{\{x,y\}\subseteq S}V(x,y)}.
$$
Define now, for $x\in \VU$, the variable  $n_x$ taking values in the set $\{0,1\}$. This variable $n_x$ can be interpreted as the occupation number of the vertex $x$: $n_x=0$ means that
the vertex is empty while $n_x=1$ means that the vertex is occupied.

\\Let $N_\VU$ be the set of all functions $\bm n:\VU\to \{0,1\}: x\mapsto n_x$ and, if $R\subset \VU$, let $N_R$ be the set of all functions $\bm n:R\to \{0,1\}$. Observe that  for each subset $R$ of $\VU$ there exists a unique function $\bm n \in N_\VU$ such that $\bm n^{-1}(1)=R$.  Therefore

$$
\begn
Z_{\GI}(z) &
=~\sum_{S\subset \VU}\,z^{|S|} ~e^{-\sum_{\{x,y\}\subset S}V(x,y)}\\
&=~\sum_{\bm n\in N_\VU}
z^{\sum_{x\in\VU}n_x}
e^{-\sum_{\{x,y\}\subset\VU}n_xn_y V(x,y)}
\egn
$$
with the convention that   $0\cdot(+\infty)=0$.

\\Expanding now the exponential, we get
\vskip.2cm
\begin{equation}
\begin{aligned}
e^{-\sum_{\{x,y\}\subset\VU}n_x n_y V(x,y)}
&= \prod_{\{x,y\}\subset\VU} [e^{- n_x n_y
V(x,y)}]\\
&= \prod_{\{x,y\}\subset\VU} [(e^{- n_x n_y
V(x,y)}-1)+1]\\
&= \sum_{k=1}^{|\VU|}\sum_{\{R_1 ,\dots ,R_k\}\in\pi(\VU)}\r(R_1)\cdots\r(R_k)
\end{aligned}
\end{equation}
where $\pi(\VU)$ is the set of all partitions of $\VU$, and
\vskip.3cm
$$
\r(R)~=~
\begin{cases} 1 &{\rm if}~ |R|=1\\
\displaystyle{\sum\limits_{g\in G_R}\prod\limits_{\{x,y\}\in
E_g}[e^{- n_x n_y V(x,y)}-1]} &{\rm if} ~ |R|\geq 2\
\end{cases}
$$
where $G_R$ is the set of connected graphs with vertex set $R$, and $E_g$ denotes the edge set of $g\in G_R$. Thus $Z_\GI(z)$ can be written as
\vv
\begin{equation}
\begin{aligned}
Z_{\GI}(z)
&=~  \sum_{\bm n\in N_\VU}
z^{\sum_{x\in\VU}n_x}\sum_{k=1}^{|\VU|}
\sum_{R_1 ,\dots,R_k\in\pi(\VU)}\r(R_1)\cdots\r(R_k)\\
&=~ \sum_{k=1}^{|\VU|}\sum_{R_1 ,\dots ,R_k\in\pi(\VU)} \sum_{\bm n\in N_\VU}
\left(\r(R_1)z^{\sum\limits_{x\in R_1}n_x}\right) \cdots
\left(\r(R_k)z^{\sum\limits_{x\in R_k}n_x} \right)\\
&=~ \sum_{k=1}^{|\VU|}\sum_{R_1 ,\dots ,R_k\in\pi(\VU)}
\prod_{i=1}^{k}{\tilde \r}(R_i,z),
\end{aligned}
\end{equation}
\vv
\\where
${\tilde \r}(R,z)= 1+z,$ if $|R|=1$, and for $|R|\geq 2$,
$$
{\tilde \r}(R,z)\,~=~\, \sum_{\bm n_R\in N_R}\;\r(R)\;z^{\sum_{x\in R}\;n_x} =
\sum_{\bm n_R\in N_R} z^{\sum_{x\in R}\;n_x} \displaystyle{\sum\limits_{g\in G_R}\prod\limits_{\{x,y\}\in
E_g}[e^{- n_x n_y V(x,y)}-1]}
$$
\\Now observe that, for any $g\in G_R$ with $|R|\ge 2$, the factor
$$
\prod_{\{x,y\}\in E_g}[e^{- n_x n_y V(x,y)}-1]
$$
is different from zero only for the configuration $\bm n_R$ such that $n_x =1$ for all $x\in R$. Therefore
$$
{\tilde \r}(R,z)~=~
\begin{cases}
1+z &{\rm if}~ |R|=1\\
z^{|R|}\sum\limits_{g\in
G_R}\prod\limits_{\{x,y\}\in E_g}[e^{-V(x,y)}-1] &{\rm if}~ |R|\geq 2
\end{cases}
$$



\\Defining  now
$$
\z(R,z)~=~
\begin{cases}
1 &{\text if}~ |R|~=~1\\
{\left(z\over 1+z\right)^{|R|}}
\sum\limits_{g\in G_R}\prod\limits_{\{x,y\}\in E_g}[e^{-V(x,y)}-1] &{\text if}~ |R|\geq 2
\end{cases}
$$
we obtain
\begin{equation}\label{eq:1zpf}
\begin{aligned}
Z_{\GI}(z) & =  (1+z)^{|\VU|}\sum_{k=1}^{|\VU|} \sum_{\{R_1 ,\dots ,R_k\}\in\pi(\VU)} \z(R_1,z)\cdots\z(R_k,z)\\
&=  (1+z)^{|\VU|}\sum_{k\ge 0}\sum_{R_1 ,\dots ,R_k \subset \VU\atop |R_i|\ge 2,~R_i\cap R_j=\0} \z(R_1,z)\cdots\z(R_k,z)
\end{aligned}
\end{equation}
where the term $k=0$ in the last sum is equal to 1 and corresponds to the partition of $\VU$ in $|\VU|$ subsets each of cardinality 1.
Finally, when  $R\subset \VU$ is such that $|R|\ge 2$, by   \reff{eq:V},  we have that
$$
\sum\limits_{g\in G_R}\prod\limits_{\{x,y\}\in E_g}[e^{-V(x,y)}-1] =\begin{cases} \sum_{E\in \mathcal{C}_{R}} (-1)^{|E|}
&~\text{if  $R\in \mathrm{C}_\VU$}\\\\0 &~\text{if  $R\notin \mathrm{C}_\VU$}
\end{cases}
$$
where we recall that $\mathcal{C}_{R}$ is the set of all connected subsets of $\EE|_R$. Therefore we have that
\be\label{eq:zxi}
Z_{\GI}(z)=(1+z)^{|\VU|}\Xi_\GI(z)
\ee
with
$$
\Xi_\GI(z)= \sum_{k\ge 0}\sum_{\{R_1 ,\dots ,R_k\}\subset \VU\atop R_i\in\mathrm{C}_\VU,~R_i\cap R_j=\0}
\z(R_1, z)\cdots\z(R_k,z)
$$
where, for any $R\in \mathrm{C}_\VU$,
$$
\z(R, z)={\left(z\over 1+z\right)^{|R|}}\sum_{E\in \mathcal{C}_{R}} (-1)^{|E|}
$$

\\The function $\Xi_\GI(z)$ can be seen as  the partition function of a polymer  gas in which set of polymers is
$\mathrm{C}_\VU$ (i.e. polymers  are connected subsets of $\VU$  with cardinality at least $2$),  each $R\in \mathrm{C}_\VU$ has
activity $\z(R, z)$, and two polymers $R, R'$ are incompatible if and only if $R\cap R'\neq \0$.

\\Let us now go back to Definition \ref{defia} and remind that given a partition scheme $\mi$ in $\GI$, $\mathcal{T}_{R}^{\mi}$ represents the set of all trees $\t$ such that $\mi(\t)=\t$. The following result is a special case of  the so-called Penrose Identity \cite{pen67}, { a central combinatorial tool which allows to rewrite alternating sums over connected subgraphs  as sums over spanning trees, with certain weights determined by the partition scheme} (see Proposition 5 in \cite{FP} {for its proof}).

\begin{lemma}\label{lem:penroseid}
    Given any partition scheme $\mi$ in $\GI$, we have, for any $R\in \mathrm{C}_\VU$,
    \begin{equation}
        \sum_{g\in \mathcal{C}_{R}} (-1)^{|E_g|} =(-1)^{|R|-1} \sum_{\t\in \mathcal{\T}^\mi_{R}} 1.
    \end{equation}
\end{lemma}

\\Hence, for a fixed partition scheme $\mi$, by Lemma \ref{lem:penroseid}, the activity of each polymer $R$, with $|R| >1$, can be rewritten as
\begin{equation}\label{eq: act}
\z(R, z)= \left(\frac{z}{z+1}\right)^{|R|}(-1)^{|R|-1}\sum_{\t\in \mathcal{\T}^\mi_{R}} 1.
\end{equation}
Therefore the sum in (\ref{eq:1zpf}) can be rewritten as a sum over forest with $k$ non-trivial trees, each one of them with vertex set $R_i$, $i=1, \cdots,k $  and such that each tree $\t$ satisfies  $\mi(\t)=\t$.
Recalling that {in (\ref{def:fmi}) we defined} $\FF_{\GI,\mi}$ {as} the set of  forests with non-trivial trees $\t$ such that $\mi(\t)=\t$, {then} we have that
$$
\Xi_\GI(z)= \sum_{F\in \FF_{\GI,\mi}}(-1)^{|F|}\left({z\over 1+z}\right)^{|V_F|}
$$
and Equation (\ref{eq:zxi}) becomes identity \reff{teoz}.

~~~~~~~~~~~~~~~~~~~~~~~~~~~~~~~~~~~~~~~~~~~~~~~~~~~~~~~~~~~~~~~~~~~~~~~~~~~~~~~~~~~~~~~~~~~~~~~~~~~~~~~~~~~~~~~~~~~~~~~~~~~~~~~~~~~~~~~~~~~~~~~~~~~~~~~~~$\Box$

\subsection*{4. Acknowledgements}
\v
The authors are grateful to Per Alexandersson for helpful bibliographic suggestions.
P.M.S.F. was supported by the Brazilian funding agency
FAPEMIG (Funda\c{c}\~ao de Amparo \`a Pesquisa do Estado de Minas Gerais).
\vv\v
\\{\bf Conflict of interest statement}. The authors declare none.
\vv\v
%
\\{\bf Data Availability Statement}.
No data was used for the research described in the article.
\vv\v
\\{\bf Funding statement}.
This work was partially supported by  Funda\c{c}\~ao de Amparo \`a Pesquisa do Estado de Minas Gerais (FAPEMIG)
under project 30105. The funder had no role in study design and analysis, decision to publish, or preparation of the manuscript.

\end{document}